\documentclass[12pt,reqno]{amsart}
\usepackage{amscd,amssymb,euscript,amstext , amsthm}
\title{The Baum-Connes Conjecture, Noncommutative Poincar\'e Duality and the Boundary
  of the Free Group}
\author{Heath Emerson}
\date{September 10, 2001}
\parskip=\baselineskip
\newcounter{axiom}[axiom]

\theoremstyle{plain}
\newtheorem{thm}{Theorem}
\newtheorem{lemma}[thm]{Lemma}
\newtheorem{prop}[thm]{Proposition}
\newtheorem{cor}[thm]{Corollary}
\newtheorem{defn}[thm]{Definition}

\newtheorem{rmk}[thm]{Remark}

\newtheorem{note}[thm]{Note}

\newcommand{\C}{\mathbb C}

\newcommand{\Z}{\mathbb Z}
\newcommand{\R}{\mathbb R}

\newcommand{\F}{\mathbb F}
\newcommand{\E}{\mathcal{E}}
\newcommand{\K}{\mathcal{ K}}

\newcommand{\Dudelta}{\hat{\Delta}}

\newcommand{\B}{\mathcal{ B}}
\newcommand{\bgamma}{\partial \Gamma}
\newcommand{\A}{C(\partial \F_{2}) \rtimes \F_{2}}
\newcommand{\BB}{C_{0}(\partial^{2} \F_{2}) \rtimes \F_{2}}

\newcommand{\Q}{\mathcal{Q}}
\begin{document}

\maketitle

\begin{abstract}

For every hyperbolic group $\Gamma$ with Gromov boundary $\bgamma$,
one can form the cross product $C^{*}$-algebra $C(\bgamma) \rtimes
\Gamma$. For each such algebra we construct a canonical
$K$-homology class, which induces a Poincar\'e duality map $K_{*}(C(\bgamma)
\rtimes \Gamma) \to K^{*+1}(C(\bgamma) \rtimes \Gamma)$. We show that
this map is an isomorphism in the case of $\Gamma = \F_{2}$ the free group on two
generators. We point out a direct connection between  our
constructions and the Baum-Connes Conjecture and eventually use the
latter to deduce our result.

\vspace{0.5cm}

2000 Mathematics Subject Classification 46L80

\end{abstract}

\section{Introduction}

The aim of this note is to point out a connection between the
Baum-Connes conjecture with coefficients for the free group $\F_{2}$
on two generators, and a Poincar\'e duality result for the `noncommutative space'
$\partial \F_{2} /\F_{2}$, where $\partial \F_{2}$ is the Gromov boundary of
$\F_{2}$, acted upon minimally by $\F_{2}$ through homeomorphisms.

In order to formulate what Poincar\'e duality should mean for a
noncommutative space such as $\partial\F_{2}/\F_{2}$, one passes to
the $C^{*}$-algebra cross product $C(\partial \F_{2}) \rtimes \F_{2}$
and to $K$-theory and $K$-homology for $C^{*}$-algebras. Poincar\'e
duality for $\partial \F_{2} / \F_{2}$ then means an isomorphism between the
$K$-theory and $K$-homology of $C(\partial \F_{2}) \rtimes \F_{2}$,
induced by cap product with a fixed $K$-homology class.

More generally one can speak of $C^{*}$-algebras having Poincar\'e
duality, or, as we call them in this paper, Poincar\'e duality algebras. It seems that such algebras are in some sense 
noncommutative analogs of spin$^{c}$ manifolds. For the commutative examples of such $C^{*}$-algebras are
given precisely by the $C^{*}$-algebras $C(M)$, where $M$ is a compact
spin$^{c}$ manifold. Such a manifold has, corresponding to the
spin$^{c}$-structure, a canonical elliptic operator on it - the Dirac
operator - and thus (see e.g. \cite{Ka1}) a canonical $K$-homology
class. Cap product with this class induces the Poincar\'e duality
isomorphism.

Various noncommutative examples of Poincar\'e duality $C^{*}$-algebras have been produced by
A. Connes, the first of which was the irrational rotation algebra
$A_{\theta}$. Several other examples now exist, but all have the same character insofar as they are
in some sense deformations of actual spin$^{c}$-manifolds. Our example
is somewhat different. The geometric data underlying $\partial \F_{2} /
\F_{2}$ is highly singular: the space $\partial \F_{2}$ is not a homology
manifold, and the group $\F_{2}$ is not a Poincar\'e duality group. It
turns out to be true, however, that in factoring the space by the action
of the group, i.e. by forming the cross product $C^{*}$-algebra
$C(\partial \F_{2}) \rtimes \F_{2}$, the resulting noncommutative space satisfies Poincar\'e duality.

Part of our goal is thus to point out this example and also to place
it in its proper context: that of hyperbolic groups acting on their
Gromov boundaries. The second part is
to show as mentioned above, a connection between our constructions
and the Baum-Connes conjecture for $\F_{2}$.

We begin by constructing - in the full generality of hyperbolic groups - the $K$-homology class cap product with which
will induce our Poincar\'e duality isomorphism. It turns out
that with Gromov hyperbolic groups $\Gamma$ in general there is a certain duality between
functions continuous on the Gromov boundary $\bgamma$ of $\Gamma$, and right translation
operators on $l^{2}\Gamma$. Using this duality,
we produce an algebra homomorphism $C(\bgamma) \rtimes \Gamma \otimes
C(\bgamma) \rtimes \Gamma \to \mathcal{Q}(l^{2}\Gamma)$,
where $\mathcal{Q}(l^{2}\Gamma) = \mathcal{B}(l^{2}\Gamma)/\mathcal{K}(l^{2}\Gamma)$ denotes the
Calkin algebra of $l^{2}\Gamma$, and where $\Gamma$ is an arbitrary
hyperbolic group. Since
$C(\bgamma) \rtimes \Gamma$ is nuclear (\cite{Del}), such an algebra
homomorphism yields via the Stinespring construction 
a class $\Delta \in KK^{1}(C(\bgamma) \rtimes \Gamma \otimes
C(\bgamma) \rtimes \Gamma , \C)$, i.e. a class $\Delta$ in the $K$-homology of
$C(\bgamma) \rtimes \Gamma \otimes C(\bgamma) \rtimes \Gamma$. Kasparov product with
$\Delta$ gives the required `cap-product' map $ \cap \Delta: K_{*}(C(\bgamma) \rtimes \Gamma) \to K^{* + 1}(C(\bgamma) \rtimes \Gamma)$.

We next wish to prove that cap product with $\Delta$ as
above gives an isomorphism in the case of $\Gamma = \F_{2}$, the
general case of hyperbolic groups being beyond the scope of this
paper. To this end  we  observe that a sort of geodesic flow
on the Cayley graph of $\F_{2}$ may be used to construct a dual
element to $\Delta$, this time    in the $K$-theory of $C(\partial \F_{2}) \rtimes \F_{2} \otimes C(\partial \F_{2}) \rtimes \F_{2}$, 
 playing the same role in this context as does the Thom class
of the normal bundle of $M$ in $M \times M$ in the commutative
setting. We obtain a putative inverse map $K^{*}(C(\partial \F_{2}) \rtimes \F_{2}) \to K_{*+1}(C(\partial \F_{2}) \rtimes \F_{2})$.

We then set about calculating the composition of these two
maps. The connection with the Baum-Connes conjecture appears in that 
the composition $K_{*}(C(\partial \F_{2}) \rtimes \F_{2}) \to K_{*}(C(\partial \F_{2}) \rtimes \F_{2})$ turns out to be
multiplication by the $\gamma$-element constructed by Julg and
Vallette. 

As mentioned, the construction of our fundamental class $\Delta$ makes
sense for a general hyperbolic group acting on its boundary, and in
fact several of our other
constructions have their counterparts for arbitrary hyperbolic
groups; thus for instance it is possible by means of work of Gromov to make sense
of `geodesic flow' for an arbitrary hyperbolic group. Furthermore,
although the statement `$\gamma = 1$' for general hyperbolic groups is
false due to the possible presence of Property T, it is nevertheless
true by work of Tu (\cite{Tu}) that $\gamma_{\bgamma \rtimes \Gamma} =
1_{C(\bgamma)}$, where $\gamma_{\bgamma \rtimes \Gamma}$ is the
  $\gamma$-element for the amenable $groupoid$ $\bgamma \rtimes \Gamma$, which weaker statement is all we need. Nevertheless,
the arguments for the general case, being substantially more involved,
will be dealt with in a later paper. We have chosen to emphasise the
free group case for two reasons: one, that the relationship to the
Baum-Connes conjecture is extremely explicit, and two, that the geometry of our
constructions is particularly visible. 

 Finally, we note that our arguments tend to suggest that the
 phenomenom of Poincar\'e duality for amenable groupoid algebras
 constructed from boundary actions of discrete groups is relatively common. Specifically, the author expects similar results
 for uniform lattices in semisimple lie groups acting on their
 Furstenberg boundaries, and for discrete, cocompact isometry groups
 of affine buildings acting on the boundaries of these
 buildings. Along these lines, we draw the attention of the reader to
 the work of Kaminker and Putnam on Cuntz-Krieger algebras  (see \cite{KP}); indeed, our result (in
 the case of the free group of two generators) can be deduced from
 theirs. In fact, our work was partly motivated by the idea of finding a
 geometric explanation for theirs.

\section{Geometric Preliminaries}

In this section we work in the generality of a Gromov hyperbolic group
$\Gamma$ (see \cite{GH} or \cite{Gro}). So let $\Gamma$ be such. Thus,
we have
fixed a generating set $S$ for $\Gamma$ and
the corresponding metric $d(\gamma_{1} , \gamma_{2}) =
|\gamma^{-1}_{1} \gamma_{2}|$, where $| \;\cdot \;  |$ denotes the word length
of a group element with respect to $S$, and with this metric $\Gamma$ is hyperbolic in
the sense of Gromov as a metric space. Note that the metric is clearly invariant under left
translation. 

Recall that with the hypothesis of hyperbolicity, the group $\Gamma$
viewed as a metric space can be compactified by addition
of a boundary: thus there exists a compact metrizable space $\bar{\Gamma} = \Gamma \cup \bgamma$
such that $\Gamma$ sits densely in $\bar{\Gamma}$, and $\bar{\Gamma}$
is compact. The compactification is $\Gamma$-equivariant in the sense
that the left translation action of $\Gamma$ extends to an action by
homeomorphisms on $\bar{\Gamma}$. 

There turns out to be an interesting duality between functions on $\Gamma$ which
extend continuously to the
Gromov compactification $\bar{\Gamma}$, and a certain class of
operators on $l^{2}\Gamma$, as follows. First we recall a
definition. For what follows, let $e_{x}$, $e_{y}$, etc,
denote the standard basis vectors in $l^{2}\Gamma$  corresponding to
points $x, y \in \Gamma$. Also, if $\tilde{f}$ is a function on
$\Gamma$, we shall denote by
$M_{\tilde{f}}$ the corresponding multiplication operator on
$l^{2}\Gamma$.

\begin{defn}\rm

An operator $T \in \mathcal{B}(l^{2}\Gamma)$ 
is finite propagation if there exists $R > 0$ such that
$<T(e_{x}) , e_{y}> = 0$ whenever $d(x,y) \ge R$. 

\end{defn}

 The duality we have alluded to is stated in the following:

\begin{lemma}

If $\tilde{f}$ is a function on $\Gamma$ which extends continuously to
$\bar{\Gamma}$, then $[M_{\tilde{f}} , T]$ is a compact operator for
all finite propagation operators $T$ on $l^{2}\Gamma$. 

\end{lemma}

For the proof, we shall need to use the following fact about the
Gromov compactification of a hyperbolic group (see \cite{GH}). 

Note that here and elsewhere in this paper, $B_{r}(x)$, for $r>0$ and
$x \in \Gamma$, denotes the ball of word-metric radius $r$ centered at
$x$.

\begin{lemma}

If $\tilde{f}$ is a continuous function on $\bar{\Gamma}$, then for
every $R>0$, we have $$\lim_{x\to \infty} \sup \{ |f(x) - f(y)| \; | \;
y \in B_{R}(x) \} = 0.$$ 

\end{lemma}

\begin{proof} (of Lemma 2)

Let $T$ be a finite propagation operator with propagation $R$, and $\tilde{f}$
a bounded function on $\Gamma$ which extends continuously to $\bar{\Gamma}$. Then $[M_{\tilde{f}} , T](e_{x}) = \sum_{y \in B_{R}(x)} \bigl(\tilde{f}(x) -
\tilde{f}(y)\bigr) T_{xy}e_{y}$ where $T_{xy}$ denotes as usual $<T(e_{x}) ,
e_{y}>$. Therefore $<[M_{\tilde{f}} , T](e_{x}) , e_{y}> = 0$ if $d(x,y) \ge
R$, and equals $\bigl(\tilde{f}(x) -
\tilde{f}(y)\bigr) T_{xy}$ else. The result follows immediately from
Lemma 3.

\end{proof}

Let $\gamma \in \Gamma$, and $\rho (\gamma)$ denote the unitary $l^{2}\Gamma \to
l^{2}\Gamma$ induced from right translation by $\gamma$, $\rho
(\gamma) e_{x} = e_{x\gamma^{-1}}$. The relevance of the above remarks
to us lies in the following observation:

\begin{lemma}

$\rho (\gamma)$ is a finite propagation operator on $l^{2}\Gamma$ for all
$\gamma \in \Gamma$. 

\end{lemma}

\begin{proof}

One has $d(x, x\gamma^{-1}) \le |\gamma|$, from which the result
follows with $R = |\gamma|$. 

\end{proof}

\begin{cor}

If $\gamma \in \Gamma$ and $\tilde{f}$ is a function on $\Gamma$ which
extends continuously to
$\bar{\Gamma}$, then $[\rho (\gamma) , M_{\tilde{f}}]$ is a compact
operator.

\end{cor}

Now, consider the unitary $I: l^{2}\Gamma \to l^{2}\Gamma$, induced from inversion $\iota:
\Gamma \to \Gamma$. Then $I \rho (\gamma) I = \lambda
(\gamma)$, where $\lambda (\gamma)$ denotes left translation by
$\gamma$; and $IM_{\tilde{f}}I = M_{\tilde{f} \circ \iota}$. The
following follows from conjugating the equation appearing in Corollary
5 by the unitary $I$:

\begin{cor}

The commutator $[\lambda (\gamma) , M_{\tilde{f} \circ
  \iota}]$ is a compact operator, for every $\gamma \in \Gamma$ and
  $\tilde{f}$ a function on $\Gamma$ extending continuously to
  $\bar{\Gamma}$.

\end{cor}

In Section $3$ we will show how the above constructions can be
organized to produce a $K$-homology class inducing a Poincar\'e
duality isomorphism.

\section{KK-theoretic preliminaries}

In this section we recall some basic facts from $KK$-theory. For
further details we refer the reader to \cite{Bla}, or to
\cite{Ka1}. 

\smallskip

{\bf KK}. $KK$ can be understood
categorically (\cite{Hig4}): there is a category
$\mathbf{KK}$ whose objects are separable, nuclear $C^{*}$-algebras and whose morphisms
$A \to B$ are the elements of $KK(A,B)$. There is a functor from
the category of $C^{*}$-algebras to the category $\mathbf{KK}$. If
$\phi:A \to B$ is an algebra homomorphism $A \to B$, we denote its image under this
functor as $[\phi]$. There
is a composition, or intersection product operation $KK(A , D)
\times KK(D , B) \to KK(A,B)$ which we denote by $(\alpha , \beta)
\mapsto \alpha \otimes_{D} \beta$. If $\phi: A \to B$ is an algebra
homomorphism, and $D$ is any $C^{*}$-algebra, we thus have a map
$\phi_{*}: KK(D, A) \to KK(D,B)$, given by $\alpha \mapsto \alpha
\otimes_{A} [\phi]$. Similarly we have a map $\phi^{*}: KK(B,D) \to
KK(A,D)$ given by $\beta \mapsto [\phi] \otimes_{B} \beta$. 

We will sometimes use the notations $\phi^{*} ([\beta])$ and $[\phi]
\otimes_{B}\beta$ interchangeably, as is warranted by clarity of
notation. Similarly with $\phi_{*}$.

If $D$ is a $C^{*}$-algebra, there is a natural map $KK(A,B) \to
KK(A \otimes D , B \otimes D)$, $\alpha \mapsto \alpha
\otimes 1_{D}$, and similarly a map $KK(A,B) \to KK(D \otimes A , D
\otimes B)$.

\smallskip

{\bf Graded Commutativity}. There are higher $KK$ groups $KK^{i}(A,B)$ for all $i \in
\Z$, defined by $KK^{i}(A , B) = KK(A , B \otimes C_{i})$ where
$C_{i}$ is the $i$th complex Clifford algebra, and one of the features of the
theory is that the intersection product is graded commutative. If $A_{1},  \ldots ,  A_{n}$ are
$C^{*}$-algebras, let $\sigma_{ij}$
denote the map $$A_{1}\otimes \cdots  A_{i} \otimes
\cdots A_{j} \otimes \cdots \otimes A_{n} \rightarrow
A_{1}\otimes \cdots  A_{j} \otimes \cdots A_{i}
\otimes \cdots \otimes A_{n}$$ obtained by flipping the
two factors. Then by graded commutativity we mean the following:
if $\alpha \in KK^{i}(A_{1},B_{1})$ and $\beta \in KK^{j}(A_{2} ,
B_{2})$, then $$(\alpha \otimes 1_{A_{2}}) \otimes_{B_{1} \otimes
  A_{2}} (1_{B_{1}}\otimes \beta )= (-1)^{ij} \,
(\sigma_{12})_{*} \sigma_{12}^{*}\bigl( ( \beta \otimes 1_{A_{1}}) \otimes
(1_{B_{2}} \otimes \alpha) \bigr) \in KK(A_{1}
\otimes  A_{2} , B_{1} \otimes B_{2}).$$

\smallskip

{\bf K-theory and K-homology}. For any $C^{*}$-algebra $A$, $KK^{i}(\C
, A) = K_{i}(A)$
is the toplogical $K$-theory of $A$, and $KK^{i}(A,\C) = K^{i}(A)$ 
is the $K$-homology of $A$ by definition.

\smallskip

{\bf Description of Even Cycles}. We let $\B(\mathcal{E})$ denote bounded
operators on a Hilbert
module $\mathcal{E}$, $\K(\mathcal{E})$ compact operators, and $\Q(\mathcal{E})$ the Calkin algebra
$\B(\mathcal{E}) / \K(\mathcal{E})$. The quotient map $\B(\mathcal{E})
\to \Q(\mathcal{E})$ will always be
denoted by $\pi$. 

Following Kasparov (\cite{Ka1}), if $\mathcal{E}$ is a Hilbert
$B$-module and $A$ acts on $\mathcal{E}$ by a homomorphism $A \to
\B(\mathcal{E})$, we will refer to $\mathcal{E}$ as a Hilbert $(A,B)$-bimodule.

Because all the algebras in this paper are
ungraded -- or alternatively, have trivial grading --  we can make certain simplifications in the
definitions of the $KK$ groups (see \cite{Bla}). With such ungraded
$A$ and $B$, cycles for $KK(A,B)$ are given simply by pairs
$(\mathcal{E} , F)$ where $\mathcal{E}$ is an $(A,B)$-bimodule, $F$ commutes modulo compact operators with the action of
$A$, and $a(F^{*}F - 1)$ and $a(FF^{*}-1)$ are compact for every $a
\in A$.

{\bf Description of Odd Cycles}. Cycles for $KK^{1}(A,B)$ are given  by pairs $(\mathcal{E} , P)$ for
which $P$ is an operator on the $(A,B)$-bimodule
$\mathcal{E}$ satisfying the three conditions
$[a,P]$, $a(P^{2} - P)$, and $a(P-P^{*})$
are compact for all $a \in A$.

Let $(\E , P)$ be an odd cycle. Then we obtain a homomorphism $A \to
\Q (\E)$ by the formula $a \mapsto \pi (PaP)$.

Conversely, let $\tau:
A \to \Q (\E)$ be a homomorphism. Under the assumption of nuclearity
of all algebras concerned, there exists a potentially larger Hilbert
$B$-module $\tilde{\E}$, a representation of $A$ on $\tilde{\E}$, an isometry $U: \mathcal{E} \to
\tilde{\mathcal{E}}$, and an operator $P$ on $\tilde{\mathcal{E}}$
such that $a(P^{2} - P)$, $[a,P]$, and $a(P-P^{*})$ are compact for all $a \in A$,
and $\pi (U^{*}PaPU) = \tau (a)$ for all $a \in A$ (see
\cite{Bla}). The data $(\tilde{\E} , P)$ makes up an odd cycle. The
process of constructing a $\tilde{\E}$, $U$, and $P$, from an
extension, we shall refer to as the Stinespring construction.

As a consequence, for $A$ and $B$ nuclear, we may regard $KK^{1}(A,B)$
as given by classes of maps $\tau: A \to \Q (\E)$, where $\E$ is a
right Hilbert $B$-module. This description of $KK^{1}$-classes will be
particularly appropriate to our purposes.

\smallskip

{\bf Bott Periodicity}. Recall that $KK^{-1}(\C , C_{0}(\R)) \cong \Z$ and is generated
by the class $[\hat{d}_{\R}]$ of the multiplier $f(x) =
\frac{x}{\sqrt{1+x^{2}}}$ of $C_{0}(\R)$, suitably interpreted in terms of
the Clifford gradings. The class
$[\hat{d}_{\R}]$ allows us to identify, for any $C^{*}$-algebras $A$
and $B$, the groups $KK^{1}(C_{0}(\R) \otimes A , B)$, and
$KK(A,B)$, by the map $KK^{1}(C_{0}(\R) \otimes A , B) \to KK(A,B)$,
$x \mapsto [\hat{d}_{\R}]\otimes_{C^{*}(\R)}
x$. We shall need to compute this map at the level of cycles in
several simple cases.

 Let $\psi$ be the function  $\psi (t) = \frac{-2i}{t+i}$ in
 $C_{0}(\R)$  It has the property that $ \psi + 1$
is unitary in $C_{0}(\R)^{+}$. We begin by stating the simplest
version of what we shall need.

\begin{lemma}

Let $\tau$ be a homomorphism $C_{0}(\R) \to \Q(H)$ to the Calkin
algebra of some Hilbert space $H$. Let
$[\tau]$ denote the class in $KK^{1}(C_{0}(\R) , \C)$ corresponding to
$\tau$. Then the class $[\hat{d}_{\R}] \otimes_{C_{0}(\R)} [\tau]
\in KK(\C , \C)$ is represented by the
cycle $(H , U+1)$, where $U$ is any operator on $H$ such that $\pi (U)
= \tau (\psi )$. 

\end{lemma}

The significance of this simple lemma is that in the given setting it is not necessary to {\it explicitly}
represent $[\tau]$ as a $KK$-cycle (that is, perform the Stinespring
construction) in order to calculate the Kasparov
product of $[\hat{d}_{\R}]$ and $[\tau]$. This is true also of the
situation in the following lemma, which will be of direct use to us.

\begin{lemma}

 Let $A_{1} , A_{2}$ be $C^{*}$-algebras and $\mathcal{E}$
 be a right Hilbert $A_{2}$-module. Let $h$ be a homomorphism
 $C_{0}(\R) \otimes A_{1} \to Q(\mathcal{E})$ and $[h]$ its class,
 regarded as an element of  $KK^{1}(C_{0}(\R)
\otimes A_{1} , A_{2})$. Assume that $h$ has the form $x \otimes a_{1} \mapsto
h'(x) h''(a_{1})$, where $h'$ and $h''$ are
homomorphisms. Suppose that the homomorphism $h''$ lifts
to a homomorphism $\tilde{h}'': A_{1} \to \B(\mathcal{E})$. Then
the class $[\hat{d}_{\R}] \otimes_{C_{0}(\R)} [h] \in KK(A_{1},A_{2})$ is
represented by the
following cycle. The module is $\mathcal{E}$ with its
 original right $A_{2}$-module structure and the left $A_{1}$-module
structure given by the homomorphism $\tilde{h}''$. The operator is
 given by $U+1$ where $U$ is any operator on
 $\mathcal{E}$ such that $\pi (U) =
 h'(\psi)$.

\end{lemma}

The proof of both lemmas involves an application of the axioms for the
intersection product, and is omitted (see \cite{Ka1}). 

\smallskip

{\bf Equivariant KK}. If $\Gamma$ is a group acting on $C^{*}$-algebras
$A$ and $B$, we have in addition to the group $KK(A,B)$, an
equivariant group $KK_{\Gamma}(A,B)$. We shall discuss this group
briefly in connection with the $\gamma$-element and the work of Julg
and Valette. Suffice it to say that the cycles for $KK_{\Gamma}(A,B)$
consist of the same cycles as for $KK(A,B)$, but with the following
extra conditions. (1) $\Gamma$ acts as linear isometric maps on the
Hilbert $(A,B)$-module $\E$, in such a way that $\gamma (a\xi b) =
\gamma(a)\gamma(\xi)\gamma(b)$ for $a \in A$, $b \in B$ and $\xi \in
\E$; (2) the operator $F$
satisfies: $\gamma (F) - F$ is compact, for all $\gamma \in \Gamma$. 

Regarding $\bf KK_{\Gamma}$ as a category in its own right, with
morphisms $A \to B$ the elements of $KK_{\Gamma}(A,B)$, and objects
$\Gamma$-$C^{*}$-algebras, there is a functor $\lambda:
KK_{\Gamma}(A,B) \to KK(A\rtimes \Gamma , B\rtimes \Gamma)$, called
{\it descent}. The map $\lambda: KK_{\Gamma}(A,B) \to KK(A\rtimes
\Gamma , B \rtimes \Gamma)$ can be explicitly calculated on cycles;
the formulas are given in \cite{Ka1}. Since $\lambda$ is a functor, it
takes the unit $1_{A} \in KK_{\Gamma}(A,A)$ to the unit $1_{A \rtimes
  \Gamma} \in KK(A\rtimes \Gamma , A\rtimes \Gamma)$, which fact we
will make use of.

\section{construction of the fundamental class}

For this section, we
shall return to the generality of a hyperbolic group $\Gamma$. Since $\Gamma$ acts by homeomorphisms on $\bgamma$, we can consider
the cross product $C^{*}$-algebra $C(\bgamma) \rtimes \Gamma$, which
is our main object of interest in this paper. Note the cross product
we are referring to is the {\it reduced} cross product; however, by
the proof of the following lemma (whose proof can be found in \cite{Del}),
the reduced and max cross products are in fact the same.

\begin{lemma}

The algebra $C(\bgamma) \rtimes \Gamma$ is nuclear and separable.

\end{lemma}

Our goal is to construct an element of the $K$-homology of the algebra
$C(\bgamma) \rtimes \Gamma \otimes C(\bgamma) \rtimes \Gamma$,
specifically, an element of $KK^{1}(C(\bgamma) \rtimes \Gamma \otimes
C(\bgamma) \rtimes \Gamma ,
\C)$. This element will be presented as an extension; that is, as a
map $C(\bgamma) \rtimes \Gamma \otimes C(\bgamma) \rtimes \Gamma \to
\Q (H)$ for some Hilbert space $H$. By our remarks in the previous
section and Lemma 9, such a map does produce a canonical class in
$KK^{1}(C(\bgamma) \rtimes \Gamma \otimes C(\bgamma) \rtimes \Gamma ,
\C)$.

We construct two commuting maps
$\lambda , \rho : C(\bgamma) \rtimes \Gamma \to \mathcal{Q}(l^{2}\Gamma)$. Let $f \in
C(\bgamma)$ and let $\tilde{f}$ denote any extension of $f$ to a
continuous function
 on $\bar{\Gamma}$. Let $M_{\tilde{f}}$ denote as
above the multiplication operator on $l^{2}\Gamma$ corresponding to $\tilde{f}$,
and let $\lambda (f)$ be the image in $\mathcal{Q}(l^{2}\Gamma)$ of
the operator
$M_{\tilde{f}}$. Let $\lambda (\gamma)$ be the image in
$\mathcal{Q}(l^{2}\Gamma)$ of the unitary $u_{\gamma}$ corresponding to left translation
by $\gamma$: $u_{\gamma}(e_{x}) = e_{\gamma x}$, $x \in \Gamma$. It is easy to check that the assignments $f \to \lambda (f)$, $\gamma \to \lambda (\gamma)$, define
a covariant pair for the $C^{*}$-dynamical system $\bigl( C(\bgamma) ,
\Gamma\bigr)$, and so a homomorphism $$\lambda: C(\bgamma) \rtimes \Gamma \to \mathcal{Q}(l^{2}\Gamma).$$ 

Next, define $$\rho: C(\bgamma) \rtimes \Gamma \to
\mathcal{Q}(l^{2}\Gamma)$$ by $\rho (a) = I \lambda (a) I$, where $I$
is at the end of Section 2.  Thus $\rho (f)$ is the image in $\mathcal{Q}(l^{2}\Gamma)$ of the
multiplication operator $M_{\tilde{f} \circ \iota}$, and $\rho
(\gamma)$ is the image in $\mathcal{Q}(l^{2}\Gamma)$ of right
translation by $\gamma$, $e_{x} \mapsto e_{x\gamma^{-1}}$. The following follows from Corollaries 5 and
6.

\begin{thm}

The homomorphisms, $\lambda , \rho: C(\bgamma) \rtimes \Gamma \to \mathcal{Q}(l^{2}\Gamma)$ commute, and so define
a homomorphism $C(\bgamma) \rtimes \Gamma \otimes C(\bgamma) \rtimes \Gamma \to \mathcal{Q}(l^{2}\Gamma)$ by $a\otimes b \to \lambda (a)
\rho (b)$. 

\end{thm}

\begin{defn}

\rm

Let $\Delta \in KK^{1}(C(\bgamma) \rtimes \Gamma \otimes C(\bgamma) \rtimes \Gamma , \C)$ denote the class
corresponding to the
above homomorphism $C(\bgamma) \rtimes \Gamma \otimes C(\bgamma) \rtimes \Gamma \to \mathcal{Q}(l^{2}\Gamma)$.

\end{defn}
 
We shall refer to the class $\Delta$ as the {\it fundamental class} of
the algebra $C(\bgamma) \rtimes \Gamma$.

Before proceeding, let us note the following. Let $\sigma_{12}: C(\bgamma) \rtimes \Gamma
\otimes C(\bgamma) \rtimes \Gamma \to C(\bgamma) \rtimes \Gamma \otimes C(\bgamma) \rtimes \Gamma$ the homomorphism which interchanges
factors, and $\sigma_{12}^{*}: KK^{1}( C(\bgamma) \rtimes \Gamma
\otimes C(\bgamma) \rtimes \Gamma, \C) \to KK^{1}( C(\bgamma) \rtimes \Gamma
\otimes C(\bgamma) \rtimes \Gamma , \C)$ the corresponding
homomorphism of $KK$ groups. The following rather simple observation reflects a common
property of `fundamental classes,' i.e. those classes implementing by
cap product Poincar\'e duality isomorphisms; the author knows of no
case, either commutative or not, where the fundamental class does not have it.

\begin{lemma}

We have: $\sigma_{12}^{*}(\Delta) = \Delta$. 

\end{lemma}

\begin{proof}

For $\sigma_{12}^{*}(\Delta)$ is the class corresponding to the map $C(\bgamma) \rtimes \Gamma
\otimes C(\bgamma) \rtimes \Gamma \to \mathcal{Q}(l^{2}\Gamma)$, $a \otimes b \mapsto \rho (a) \lambda (b)$. But
this is unitarily conjugate to the map $a \otimes b \mapsto \lambda
(a) \rho (b)$ via the symmetry $I$. 

\end{proof}

We can now define the `cap-product map' interchanging the $K$-theory
and $K$-homology of $C(\bgamma) \rtimes \Gamma$, which we are going to show is an isomorphism
when $\Gamma = \F_{2}$. Specifically, define: $$\cap\Delta:
K_{*}(C(\bgamma) \rtimes \Gamma)
\to K^{*+1}(C(\bgamma) \rtimes \Gamma)$$by the formula $$x \mapsto
(x\otimes 1_{C(\bgamma) \rtimes \Gamma})
\otimes_{C(\bgamma) \rtimes \Gamma \otimes C(\bgamma) \rtimes \Gamma}\Delta.$$

Our main theorem is the following:

\begin{thm}

For $\Gamma = \F_{2}$ and $\Delta$ as in Defintion 11, the map $\cap \Delta$ is an
isomorphism.

\end{thm}


\section{Connes' notion of Poincar\'e duality}

In order to prove that the map $\cap \Delta$ of the previous section
is an isomorphism, we shall use some ideas due to Connes.

\begin{thm}

Let $A$ be a separable, nuclear $C^{*}$-algebra, and $\Delta $
a class in $KK^{i}(A \otimes A , \C)$. 

Suppose we can find a class $\Dudelta \in KK^{-i}(\C , A \otimes
A)$ such that the following equations hold:

\begin{equation} \bigl( \Dudelta \otimes 1_{A} \bigr)
  \otimes_{A \otimes A \otimes A}\bigl(1_{A} \otimes  \sigma_{12}^{*}\Delta \bigr)
  = 1_{A}, \end{equation} and
\begin{equation}\bigl(  (\sigma_{12})_{*}\Dudelta  \otimes
  1_{A}\bigr)\otimes_{A\otimes A\otimes A}\bigl(1_{A} \otimes \Delta
  \bigr) = (-1)^{i} \;  1_{A}. \end{equation}

Then the map $$\cap \Delta: K_{j}(A) \to K^{j+i}(A)$$ defined previously, is an
isomorphism with inverse (up to sign) the map $K^{j}(A) \to
K_{j-i}(A),$ $$y \mapsto \Dudelta \otimes_{A \otimes A} (1_{A} \otimes
y).$$

\end{thm}

If $A$ is as above, with classes $\Delta$ and $\Dudelta$ satisfying
Equations (1) and (2), we will call $A$ a Poincar\'e duality algebra.

\begin{proof}

The hypotheses imply the two equations: \begin{equation} \bigl( \Dudelta \otimes 1_{A} \bigr) \otimes_{A \otimes A \otimes
    A} \bigl ( 1_{A} \otimes \sigma_{12}^{*}(\Delta) \bigr) = 1_{A}, \end{equation}
    and \begin{equation} \bigl( (\sigma_{12})_{*}(\Dudelta) \otimes 1_{A} \bigr) \otimes_{A \otimes A \otimes
    A} \bigl ( 1_{A} \otimes \Delta \bigr) = (-1)^{i} \; 1_{A}. \end{equation}

We show that as a consequence of these two equations, \begin{equation}
  \Dudelta \otimes_{A
  \otimes A} \bigl(1_{A} \otimes (y  \cap \Delta )\bigr) = (-1)^{ij} \; y,
  \; y \in KK^{j}(\C , A) \end{equation}

Expanding the product involved in $(5)$, we obtain: $$\Dudelta \otimes_{A \otimes A}
\bigl( 1_{A} \otimes y \otimes 1_{A} \bigr) \otimes_{A \otimes A
  \otimes A} (1_{A} \otimes \Delta).$$ 

Consider the term $(1_{A}
\otimes y \otimes 1_{A}) \otimes_{A \otimes A \otimes A} (1_{A}
\otimes \Delta).$ It is easy to check this is the same as $\bigl(
1_{A} \otimes 1_{A} \otimes y \bigr) \otimes_{A \otimes A \otimes A}
\bigl( 1_{A} \otimes \sigma_{12}^{*}(\Delta)\bigr).$ Returning to the
original product (5), we see the latter can be written

$$\bigl(\Dudelta
\otimes_{A\otimes A} (1_{A\otimes A}\otimes y) \bigr) \otimes_{A \otimes A \otimes A} \bigl( 1_{A} \otimes
\sigma_{12}^{*}(\Delta)\bigr).$$

Now, by skew-commutativity of the
external tensor product, $$\Dudelta \otimes_{A\otimes A} (1_{A\otimes
  A}\otimes y) = (-1)^{ij}
(\sigma_{23})_{*}(\sigma_{12})_{*} \bigl( y \otimes_{A}(1_{A} \otimes
\Dudelta )
\bigr) = (-1)^{ij} y \otimes_{A} 
(\sigma_{23})_{*}(\sigma_{12})_{*} (1_{A} \otimes \Dudelta
).$$Furthermore, $(\sigma_{23})_{*}(\sigma_{12})_{*} (1_{A} \otimes \Dudelta ) =
\Dudelta \otimes 1_{A}$. Hence, putting back into the main product, we see that (5) can be written $$(-1)^{ij} \; y
\otimes_{A} \bigl( (\Dudelta \otimes 1_{A} ) \otimes_{A \otimes A
  \otimes A} (1_{A} \otimes \sigma_{12}^{*}(\Delta) \bigr) = (-1)^{ij}
\; y,$$where the last equality follows from equation $(1)$.

\end{proof}

\begin{rmk}

\rm

We note that if we happen to have $\Delta$ and $\Dudelta$ as above,
and  $$\sigma_{12}^{*}(\Delta) = \Delta$$ and
$$(\sigma_{12})_{*}(\Dudelta) = (-1)^{i}\; \Dudelta,$$then the two
equations $(1)$ and $(2)$ above would be the same, and it would suffice to
show that one of them holds. This is the case in the commutative
setting of a compact spin$^{c}$-manifold, and will be the case for us,
also, part of which we have already proven (Lemma 12).

\end{rmk}

We now set about proving Theorem 13 in the case of $\Gamma = \F_{2}$ by verifying the equations $(1)$
and $(2)$ of Theorem 14 above, with, i.e. $A = \A$ and $\Delta$ the fundamental
class of Definition 11. We need first produce an element $\Dudelta \in
KK^{-1}( \C , \A \otimes \A)$ playing the role of the dual
element in Theorem 14. We will then verify equation (1), the other being rendered superfluous as a consequence of
  Remark 15, which is applicable in this case. 

It will turn out, rather surprisingly, that equation (1) can be shown
to be equivalent to the equation  $$\gamma_{\partial \F_{2} \rtimes \F_{2}}  =
    1_{C(\partial \F_{2})},$$where $\gamma_{\partial \F_{2} \rtimes
      \F_{2}}$ is the $\gamma$-element for the groupoid $\partial \F_{2}
      \rtimes \F_{2}$. Since this latter equation has been established by Julg and
      Valette, and also by J.L. Tu, we will by this device, i.e. by
      means of
      the Baum-Connes Conjecture, be done.

\section{Construction of a dual element}

In this section as for the rest of this note we specialize to the free
group $\F_{2}$ on two generators. We are going to
define an element $\Dudelta \in KK^{-1}(\C , \A \otimes \A)$ serving
as an `inverse' to $\Delta$. 

 $\Dudelta$ shall be constructed by use of the fact that any two points of $\partial \F_{2}$
may be connected by a unique geodesic.

By ``geodesic'' we shall mean an isometric map $r: \Z \to
\F_{2}$. Topologize the collection of such $r$ by means
of the metric $$d_{G\F_{2}}(r_{1} , r_{2}) = \sum_{n\in \Z} 2^{-|n|} d(r_{1}(n) ,
r_{2}(n))$$and denote the resulting metric space by $G\F_{2}$ (we
follow \cite{Gro}). Both
$\F_{2}$ and $\Z$ act freely and properly on $G\F_{2}$, the former by
translation $(\gamma r) (n) = \gamma r(n)$, and the latter by flow
$(g^{n}r) (k) = r(k- n).$ These actions commute. Note that $G\F_{2}/\F_{2}$ is compact,
whereas $G\F_{2} /\Z$ may be identified with the $\F_{2}$-space
$$\partial^{2}\F_{2} = \{ (a,b) \in \partial \F_{2} \times \partial \F_{2} \; | \; a
\not= b\}.$$All these observations are easy to check. As a consequence
of them, the $C^{*}$-algebras $C(G\F_{2} / \F_{2}) \rtimes \Z$ and
$C_{0}(\partial ^{2}\F_{2}) \rtimes \F_{2}$ are strongly Morita
equivalent (see \cite{Rie}). Let $[E]$ denote the class of the strong
Morita equivalence bimodule. It is an element of
$KK(C(G\F_{2}/\F_{2})\rtimes \Z , C_{0}(\partial ^{2}\F_{2}) \rtimes
\F_{2})$. 

On the other hand, if $u$ is the generator of $\Z \subset C^{*}(\Z)
\subset C(G\F_{2}/\F_{2}) \rtimes \Z$, we obtain a natural
homomorphism $C_{0}(\R) \to C(G\F_{2} /\F_{2}) \rtimes \Z$ by the
formula $\psi \mapsto u-1$ where, recall, $\psi$ is a specified generator of
$C_{0}(\R)$ satisfying $\psi+1 \in C_{0}(\R)^{+}$ is unitary.

We denote the class in $KK(C_{0}(\R) , C(G\F_{2}/\F_{2})\rtimes \Z)$
of this homomorphism by $[u-1]$.

It will be convenient for our later computations to define an
auxilliary class $[D]$, which will lie in $KK^{-1}(\C ,
C_{0}(\partial^{2}\F_{2}) \rtimes \F_{2})$, as follows.

\begin{defn}

\rm

The class $[D] \in KK^{-1}(\C , C_{0}(\partial ^{2}\F_{2}) \rtimes
\F_{2})$ shall be defined by $$[D] = [\hat{d_{\R}}]
\otimes_{C_{0}(\R)}[u-1] \otimes_{C(G\F_{2}/\F_{2})\rtimes \Z} [E].$$

\end{defn}

Next, note that the cross product $C_{0}(\partial^{2}\F_{2})
\rtimes \F_{2}$ may be regarded as a subalgebra of $C(\partial \F_{2})
\rtimes \F_{2} \otimes C(\partial \F_{2})
\rtimes \F_{2}$, via
the composition of inclusions: $$C_{0}(\partial^{2}\F_{2})
\rtimes \F_{2} \to C(\partial \F_{2} \times
\partial \F_{2}) \rtimes \F_{2} \cong \bigl( C(\partial \F_{2}) \otimes C(\partial \F_{2})\bigr) \;
\rtimes \F_{2} \to C(\partial \F_{2}) \rtimes \F_{2} \otimes
C(\partial \F_{2}) \rtimes \F_{2}.$$ Let $i$ denote this composition.

Our class $\Dudelta$ will be defined by:

\begin{defn}

\rm

Let $\Dudelta = [D] \otimes_{C_{0}(\partial ^{2}\F_{2}) \rtimes
  \F_{2}}[i] \in KK^{-1}(\C , \A \otimes \A),$ where $[\hat{d}_{\R}]$
  is as in Section 3, and $[u-1]$ and $[E]$ are as above. 

\end{defn}

It will be convenient to calculate more explicitly the cycle corresponding
to the class
$[u-1]\otimes_{C(G\F_{2}/\F_{2})\rtimes \Z} [E] \in KK(C_{0}(\R) ,
C_{0}(\partial ^{2}\F_{2}) \rtimes \F_{2})$. We will express it as a
homomorphism $C_{0}(\R) \to C_{0}(\partial ^{2}\F_{2}) \rtimes
\F_{2}$; that is, as an element $w \in C_{0}(\partial ^{2}\F_{2})
\rtimes \F_{2}$ such that $w+1$ is unitary in $\bigl(C_{0}(\partial^{2}\F_{2})
\rtimes \F_{2}\bigr)^{+}$. 

We will first describe an element $v \in C_{0}(\partial
^{2}\F_{2}) \rtimes \F_{2}$ satisfying $v^{*}v = vv^{*} =
\chi$, where $\chi$ is a projection. We will then set $w =
v-\chi$. Then, of course, $w+1 = v +1-\chi$ will be unitary in
$\bigl( C_{0}(\partial
^{2}\F_{2}) \rtimes \F_{2}\bigr)^{+}$. 

 As the method of discovering such an explicit
description (that is, of transfering $K$-classes under strong Morita
equivalences) is well known (see \cite{Co2} in which a similar
calculation is carried out in the context of $A_{\theta}$) we give the outcome without
further discussion.

As a function on
$\partial^{2}\F_{2} \times \F_{2}$, $v(a,b, \gamma) = 1$ if and only if there
exists a geodesic $r_{a,b}$ such that $r_{a , b}(-\infty) = a$,
$r_{a,b}(+\infty) = b$, $r_{a,b}(0) = e$, and $r_{a,b}(-1) =
\gamma$. And $v(a,b, \gamma) = 0$ else.

Note that  $\chi = v^{*}v = vv^{*}$ is the locally constant function on $\partial^{2}\F_{2}$
  given by $\chi(a,b) = 1$ if some (therefore any) geodesic from $a$
  to $b$ passes through $e$, and equals $0$ else.

We can describe $v$ in group-algebra notation as follows. Fix $\gamma$
a generator. Then $v(\cdot \; , \cdot \; , \gamma)$ is a function on
$\partial^{2}\F_{2}$, and in particular is a function on $\partial \F_{2}
\times \partial \F_{2}$, whose representation as a tensor product of two
functions on $\partial \F_{2}$ 
is: $$v(\cdot \; , \cdot \; , \gamma) = \chi_{\gamma}
\otimes (1-\chi_{\gamma}),$$ where $$\chi_{\gamma}(a) =
\begin{cases}  1 &  \quad \gamma \in [e,a)\\0& \quad \mathrm{else} \end{cases}.$$

We can therefore represent $v$ as $$v = \sum_{|\gamma| = 1}
\chi_{\gamma} \gamma \otimes (1-\chi_{\gamma})\gamma \in
C_{0}(\partial^{2}\F_{2})\rtimes \F_{2} \subset  C(\partial
\F_{2}) \rtimes \F_{2} \otimes
  C(\partial \F_{2}) \rtimes \F_{2}.$$Similarly we we represent the
  function $\chi$ by $\chi = \sum \chi_{\gamma} \otimes
  (1-\chi_{\gamma}),$ and it is easy to check that $v^{*}v = vv^{*} = \chi$, as
  claimed.

Finally, we note the following:

\begin{lemma}

The class $\Dudelta$ satisfies $(\sigma_{12})_{*}(\Dudelta) =
-\Dudelta$.

\end{lemma}

\begin{proof}

We have $\Dudelta = i_{*}([D])$, and so $(\sigma_{12})_{*}(\Dudelta)
= (\sigma_{12})_{*}i_{*}([D]) = (\sigma_{12} \circ i )_{*}([D]) =
(\bar{\sigma}_{12})_{*}([D])$, where $\bar{\sigma}_{12}:
C_{0}(\partial^{2}\F_{2})\rtimes\F_{2} \to \BB$ is the algebra homomorphism induced
by the $\F_{2}$-equivariant map $\partial^{2}\F_{2} \to
\partial^{2}\F_{2} $, $(a,b) \mapsto (b,a)$. Now $[D] = [\hat{d}_{\R}]
\otimes_{C_{0}(\R)} [v-\chi]$ and hence $(\bar{\sigma}_{12})_{*}([D])
= [\hat{d}_{\R}] \otimes_{C_{0}(\R)} (\bar{\sigma}_{12})_{*}
([v-\chi]) = [\hat{d}_{\R}] \otimes_{C_{0}(\R)} [v^{*}-\chi] =
-[\hat{d}_{\R}] \otimes_{C_{0}(\R)} [v^{*}-\chi]$ by a direct calculation and we are done.

\end{proof}

In the following sections we will show that in an appropriate sense
$\Dudelta$ provides an `inverse' to the extension $\Delta$. More precisely,
we will show that the conditions of Theorem 14 are
met by $\Delta$ the fundamental class, and the element $\Dudelta$
above.

\section{The $\gamma$-element}

Before proceeding to verify the equations of Theorem 14, we will need
to recall the work of Julg and Vallette (\cite{JV}).

Up to now we have adopted the convention of writing even $KK$-cycles in the
form $(\E , F)$, where $F$ is an operator on the module $\E$. A
different definition is possible, in which {\it two} modules are
involved, and $F$ is an operator between them. This was the set-up in
the paper of \cite{JV}. We will retain their notation temporarily. In
a moment we will describe a means of {\it geometrically} describing
their class in a way consistent with our conventions.

 Consider the Cayley graph $\Sigma$ for
$\F_{2}$, which is a tree with edges $\Sigma^{1}$ and vertices
$\Sigma^{0}$. Note that we work with $geometric$ edges, i.e. set
theoretic pairs of vertices $\{x, x'\}$. If $x$ is a
vertex, let $x'$ be the vertex one unit closer to $e$, the origin, and
let $s(x)$ be the edge $\{ x , x'\}$. Define an operator $$b:
l^{2}\Sigma^{0} \to l^{2}\Sigma^{1}$$by $$b(e_{x}) = \begin{cases}
  e_{s(x)} &  \quad x \not= e\\0& \quad x=e \end{cases}.$$Then it is
clear that $b$ is an isometry, is Fredholm, and has index
$1$. Next, note that $\F_{2}$ acts unitarily on $l^{2}(\Sigma^{0})$
and $l^{2}(\Sigma^{1})$, and that, furthermore, $\gamma b \gamma^{-1}
- b$ is a compact (in fact finite rank) operator, for all $\gamma \in
\F_{2}$.

It follows that the pair $\bigl( l^{2}\Sigma^{0} \oplus
l^{2}\Sigma^{1} , \begin{pmatrix} 0 & b^{*} \\ b & 0\end{pmatrix}\bigr)$
defines a cycle for $KK_{\F_{2}}(\C , \C)$.

Let $\gamma$ denote its
class. That $\gamma = 1$ in this group implies the Baum-Connes
conjecture for $\F_{2}$. This fact (that $\gamma = 1$) was proved by
Julg and Valette in \cite{JV}.

We can produce a cycle for $KK_{\F_{2}}(C(\partial \F_{2}) , C(\partial \F_{2}))$,
whose class we will denote by $\gamma_{\partial \F_{2}}$, by tensoring all the
above data with $C(\partial \F_{2})$. Thus, let $\mathcal{E}^{0} = C(\partial \F_{2}
; l^{2}(\Sigma^{0}))$, and $\mathcal{E}^{1} = C(\partial \F_{2} ;l^{2}(\Sigma^{1}))$. Let
$B: \mathcal{E}^{0} \to \mathcal{E}^{1}$ be defined by $(B\xi)(a) =
b(\xi (a))$. The Hilbert $C(\partial \F_{2})$-modules $\mathcal{E}^{i}$ carry
obvious actions of $\F_{2}$. Let $\gamma_{\partial \F_{2}}$ be the class of
the cycle $\bigl( \mathcal{E}^{0} \oplus \mathcal{E}^{1} ,
\begin{pmatrix} 0& B^{*} \\ B & 0\end{pmatrix})$. It is easy to check
that the process of tensoring with $C(\partial \F_{2})$ in this way preserves
units; that is:  $$\gamma = 1 \Rightarrow \gamma_{\partial \F_{2}} = 1_{C(\partial \F_{2})}$$in the ring
$KK_{\F_{2}}(C(\partial \F_{2}) , C(\partial \F_{2}))$. Hence, we
have:

\begin{lemma}

The cycle $\bigl( \mathcal{E}^{0} \oplus \mathcal{E}^{1} ,
\begin{pmatrix} 0& B^{*} \\ B & 0\end{pmatrix})$ is equivalent to the
cycle corresponding to 
$1_{C(\partial \F_{2})}$ in the group $KK_{\F_{2}}(C(\partial \F_{2}) , C(\partial \F_{2}))$. 

\end{lemma}

We now set about describing a cycle equivalent to the above but which
is in some sense simpler. To do this it will be notationally and
conceptually simpler to work with
fields. Thus, we note that $\mathcal{E}^{0}$ and $\mathcal{E}^{1}$ may be
viewed as sections of the constant fields of Hilbert spaces
$\{H^{0}_{a} \; | \; a \in \partial \F_{2} \}$, respectively $\{H^{1}_{a} \; | \; a \in
\partial \F_{2}\}$, with
$H^{0}_{a} = l^{2}(\Sigma^{0})$ and $H^{1}_{a} =
l^{2}(\Sigma^{1})$ for all $a \in \partial \F_{2}$, and that the operator $B$
may be regarded as the constant family of operators $\{ b_{a}\; | \; a
\in \partial \F_{2} \}$ with $b_{a} = b$ for
all $a \in \partial \F_{2}$. What we are going to do is eliminate
edges from the cycle at the expense of changing the constant field of
operators to a nonconstant field.

To this end consider the field of unitary maps $\{ U_{a}:
H^{1}_{a} \to H^{0}_{a}\; | \; a \in \partial \F_{2} \}$ given by $U_{a} (e_{s})
= e_{x}$, where $x$ is the vertex of $s$ farthest from $a$. Note that the
assignment $a \mapsto U_{a}$, though not constant, is strongly continuous. For if $a$ and $b$ are two boundary
points, then $U_{a} = U_{b}$ except for edges lying on the geodesic
$(a,b)$. Consequently, if $s$ is a fixed edge, and $a$ and $b$ are
close enough, then $U_{a}(e_{s}) = U_{b}(e_{s})$, since if $a$ and $b$ are
sufficiently close, $s$ does not lie on $(a,b)$. 

Now, consider the composition $$l^{2}\F_{2} = H^{0}_{a}
\stackrel{b_{a} }{\longrightarrow} H^{1}_{a} \stackrel{U_{a}}
{\longrightarrow} H^{0}_{a} = l^{2}\F_{2},$$ which we denote by
$W_{a}$. We see that for $x = e$, $W_{a}(e_{x}) = 0$, and for $x \not= e$
we have:
$$W_{a}(e_{x}) = \begin{cases}
  e_{x'} &  \quad x\in [e,a) \\e_{x}& \quad \mathrm{else}
\end{cases},$$where $x'$ is the vertex one unit closer to $e$ than
$x$.

Since the assignment $a \to W_{a}$ is continuous, we obtain a Hilbert
$C(\partial \F_{2})$-module map $\mathcal{E}^{0} \to \mathcal{E}^{0}$ by
defining for $\xi \in C(\partial \F_{2} ;
l^{2}\F_{2})$, $(W\xi)(a) = W_{a}(\xi (a))$. Then, by unitary
invariance of $KK$ and the work of Julg and Vallette, we see:

\begin{lemma} The cycle $\bigl( \mathcal{E}^{0}
\oplus \mathcal{E}^{0} , \begin{pmatrix} 0& W^{*}\\W &
  0\end{pmatrix}\bigr)$ is equivalent to the cycle corresponding to $1_{\partial \F_{2}}$ in the group
  $KK_{\F_{2}}(C(\partial \F_{2}) , C(\partial \F_{2}))$. 

\end{lemma}

Since we have now altered the cycle of Julg and Valette up to
equivalence so that only one
Hilbert module is involved (it is now otherwise known as an `evenly graded'
Fredholm module), we may now return as promised to our conventions and write it simply
$$( C(\partial \F_{2} ; l^{2}\F_{2})    , W),$$consistent with the way we have been
writing (even) $KK$-cycles up to now.

To summarize, we have: $$[( C(\partial
\F_{2} ; l^{2}\F_{2})    , W)] = [1_{C(\partial \F_{2})}] \in
  KK_{\F_{2}}(C(\partial \F_{2}) , C(\partial \F_{2})).$$

We shall next apply the descent map $$\lambda: KK_{\F_{2}}(C(\partial \F_{2}) ,
C(\partial \F_{2})) \to KK(C(\partial \F_{2}) \rtimes
\F_{2},C(\partial \F_{2}) \rtimes \F_{2})$$ to the cycle described
above, thus producing
a cycle for $KK(C(\partial \F_{2}) \rtimes \F_{2},C(\partial \F_{2})
\rtimes \F_{2})$ which by functoriality of descent will be equivalent
to the cycle corresponding to $1_{C(\partial \F_{2})
  \rtimes \F_{2}}$. A direct appeal to the definition of $\lambda$
(see \cite{Ka1}) produces the cycle  $\bigl( C(\partial \F_{2}) \rtimes \F_{2} \otimes
l^{2}\F_{2} , \bar{W} \bigr)$, where,
    regarding $C(\partial \F_{2}) \rtimes \F_{2} \otimes l^{2}\F_{2}$ as given by 
    functions $\F_{2} \to C(\partial \F_{2}) \otimes l^{2}\F_{2}$, the action of $\bar{W}$ on these
    functions is given by the formula $(\bar{W}\xi )(\gamma) = W(\xi
    (\gamma))$. We have:

\begin{lemma}

The cycle $\bigl( C(\partial \F_{2}) \rtimes \F_{2} \otimes
l^{2}\F_{2} , \bar{W} \bigr)$ is equivalent to the cycle
corresponding to $1_{C(\partial \F_{2}) \rtimes \F_{2}}$ in $KK(C(\partial \F_{2}) \rtimes \F_{2},C(\partial \F_{2}) \rtimes \F_{2})$. 

\end{lemma}

This concludes our preparatory work. We will now show that the class
of the cycle given in the above lemma is the same as the class of the Kasparov
product of the elements $\Dudelta$ and $\Delta$, concluding thus as a
consequence of the work of Julg and Valette that equation (1) holds.

\section{untwisting}

We are interested in calculating the cycle corresponding to the
Kasparov product $$\bigl(\Dudelta \otimes 1_{C(\partial \F_{2}) \rtimes \F_{2}} \big) \otimes_{C(\partial \F_{2}) \rtimes \F_{2}\otimes C(\partial \F_{2}) \rtimes \F_{2} \otimes C(\partial \F_{2}) \rtimes \F_{2}}
\bigl( 1_{C(\partial \F_{2}) \rtimes \F_{2}} \otimes
\sigma_{12}^{*}\Delta\bigr).$$

In this section we will do something we call - following an analogous
procedure in \cite{KP} - `untwisting.' A simple but fundamental
property of hyperbolic groups - and in particular of the free group -
will be used: specifically, {\it if two points $a $ and $b$ on $\partial
\F_{2}$ are sufficiently far apart then any
geodesic connecting them passes quite close to the identity $e$ of the
group}. This follows immediately from the definition of the topology on the
compactified space $\overline{\F}_{2}$. More precisely:

\begin{lemma}

Let $\tilde{N} $ be a
neighbourhood of the diagonal $\{(a,a) \; | \; a \in \partial
\F_{2}\}$ in $ \partial \F_{2} \times \overline{\F}_{2}$. Then there exists $R>0$ such that if $(a,b) \in (\partial
\F_{2} \times \overline{\F}_{2}) \backslash \tilde{N}$, then the (unique)
geodesic from $a$ to $b$ passes through $B_{R}(e)$.

\end{lemma}

\begin{note}

\rm

To simplify notation in this section, we shall denote by $A$ the cross
product $C(\partial \F_{2}) \rtimes \F_{2}$, and by $B$ the algebra
$C_{0}(\partial^{2}\F_{2}) \rtimes \F_{2}$. 

\end{note}

Consider then the product $ \bigl(\Dudelta \otimes 1_{A} \big) \otimes_{A\otimes A \otimes A}
\bigl( 1_{A} \otimes \sigma_{12}^{*}\Delta\bigr)$ involved on the left
hand side of equation (1).

Since $\Dudelta = i_{*}([D]) = [D] \otimes_{B} [i]$, we have $$\bigl(\Dudelta \otimes 1_{A} \big) \otimes_{A\otimes A \otimes A}
\bigl( 1_{A} \otimes \sigma_{12}^{*}\Delta\bigr) = \bigl( [D] \otimes 1_{A}\bigr) \otimes_{B
    \otimes A} [i\otimes 1_{A}]\otimes_{A\otimes A \otimes A} (1_{A} \otimes
  \sigma_{12}^{*}\Delta).$$

We will begin by examining the term $[i\otimes 1_{A}]\otimes_{A\otimes
  A \otimes A} (1_{A} \otimes
  \sigma_{12}^{*}\Delta) \in KK^{1}(B
  \otimes A , A)$. It is easy to describe the corresponding cycle explicitly. For since $\sigma_{12}^{*}\Delta$ is represented by a map
  $A \otimes A \to \mathcal{Q}(l^{2}\F_{2})$, so also $1_{A} \otimes \sigma_{12}^{*}\Delta$ is
  represented by a map $A \otimes A \otimes A \to \mathcal{Q}(A \otimes l^{2}\F_{2})$, and $[i
  \otimes 1_{A}]\otimes_{A\otimes A \otimes A} (1_{A} \otimes \sigma_{12}^{*}\Delta)$ is
  represented by a map $B \otimes A \to \mathcal{Q}(A \otimes
  l^{2}\F_{2})$. By construction, this map is given on the set of elementary tensors by
  the formula  \begin{equation} a_{1} \otimes a_{2} \otimes a_{3} \mapsto a_{1} \otimes \rho
  (a_{2}) \lambda (a_{3}),\end{equation} where we have suppressed the
  inclusion $i: B
  \to A \otimes A$, so that $a_{1} \otimes a_{2}$ in the above
  expression is understood as an element of $B$. 

We shall first show that the above map up to unitary equivalence can
be rewritten in a much more tractable way.

Before proceeding, let $\tilde{G}$ be a function on $\partial \F_{2} \times
\F_{2}$ not necessarily continuous in the second variable, but
continuous in the first. Then $\tilde{G}$ can be made to act on the
right $A$-module $A \otimes l^{2}\F_{2}$ by the formula $$\tilde{G}\cdot (a \otimes e_{y}) =
\tilde{G}(\; \cdot \; 
, y)a \otimes e_{y},$$noting that for each $y \in \F_{2}$, $\tilde{G}(\cdot \; , y) \in
C(\partial \F_{2}) \subset A$. 

Now let $F$ be
a continuous, compactly supported function on
$\partial^{2}\F_{2}$. Thus $F$ is a
continuous function on $\partial \F_{2} \times \partial  \F_{2}$ vanishing in a
neighbourhood of the diagonal. So we can extend it to a continuous
function $\tilde{F}$ on $\partial \F_{2} \times \overline{\F}_{2}$ by the Tietze
Extension Theorem and restrict the result to $\partial \F_{2} \times
\F_{2}$. Let
$\tilde{F}'$ denote the function on $\partial \F_{2} \times \F_{2}$
given by $(a,x) \mapsto \tilde{F}(x^{-1}a , x^{-1})$. Note that
$\tilde{F}'$ is continuous in the first variable but not in the
second. Hence $\tilde{F}'$ may be made to act on the Hilbert
$A$-module $A \otimes l^{2}\F_{2}$ by the remark in the previous
paragraph.  We can thus regard $\tilde{F}'$ as an element of $\B (A
\otimes l^{2}\F_{2})$. Let $\tau (F)$ denote the image of the
operator $\tilde{F}'$ in $\Q (A \otimes l^{2}\F_{2})$.

Remark that $F \mapsto \tau (F)$ is a well-defined homomorphism
$C_{0}(\partial^{2}\F_{2}) \to \Q(A \otimes l^{2}\F_{2})$. For any two
extensions of $F$ to functions on $\partial \F_{2} \times
\overline{\F}_{2}$ differ by a function - say $\tilde{H}$ - vanishing on
$\partial \F_{2} \times \partial \F_{2}$. Then $\tilde{H}'$ also
vanishes on $\partial \F_{2} \times \partial \F_{2}$, and so defines
an operator lying in $\K (A \otimes l^{2}\F_{2})$. 

Next, for $\gamma \in \F_{2}$, set 
$\tau (\gamma ) =  1 \otimes \rho(\gamma) \in \mathcal{Q}(A
  \otimes l^{2}\F_{2}).$ It is a routine computation to check that the
  assignments $$F \mapsto \tau (F)$$and $$\gamma \mapsto \rho
  (\gamma)$$ make up a covariant pair for the dynamical system $\bigl(
  C_{0}(\partial^{2}\F_{2}) , \F_{2}\bigr)$, and hence a homomorphism
  $$\tau: B \to \mathcal{Q}(A \otimes l^{2}\F_{2}).$$

Next, define
a covariant pair for the dynamical system $(C(\partial \F_{2}) , \F_{2})$ by
$\varphi (f) = f \otimes 1 \in \mathcal{B}(A \otimes l^{2}\F_{2})$, and $\varphi (\gamma) = \gamma \otimes
 u_{\gamma} \in \mathcal{B}( A \otimes l^{2}\F_{2})$. It is similarly easy to check this
 makes up a covariant pair and so a homomorphism $$\varphi: A \to \mathcal{B}(A
 \otimes l^{2}\F_{2}).$$

 The following proposition is key to the untwisting argument.

\begin{prop}

The class $[i\otimes 1_{A}] \otimes_{A\otimes A \otimes A}\bigl( 1_{A} \otimes
    \sigma_{12}^{*}\Delta\bigr) \in
  KK^{1}(B\otimes A , A)$ is represented by the homomorphism $\iota: B
  \otimes A \to \mathcal{Q}(A \otimes l^{2}\F_{2})$, $$\iota (b \otimes a ) = \tau
  (b)\pi (\varphi (a)), \; \; b \in B , a \in A$$ where
    $\varphi$, $\tau$ are as above.

\end{prop}

We note that the homomorphisms $\tau$ and $\pi \circ \varphi$ commute,
and so $\iota$ actually is a homomorphism. That $\iota$ is a
homomorphism also follows, however, from the proof of Lemma 24 below,
which shows that $\iota$ is unitarily conjugate to the map in Equation
6.

We will require the following: 

\begin{lemma} Let $k \in C_{c}(\partial^{2}\F_{2} \times \partial \F_{2})$, and
  $\tilde{k}$ an extension of $k$ to a continuous function on $\partial \F_{2}
  \times \overline{\F}_{2} \times \overline{\F}_{2}$. Then the two functions on
  $\partial \F_{2} \times \F_{2}$  $$(a,x) \mapsto
  \tilde{k}(x^{-1}(a) , x^{-1} , x) $$and $$(a,x) \mapsto
  \tilde{k}(x^{-1}(a) , x^{-1} , a)$$are the same modulo $C_{0}(\partial \F_{2}
  \times \F_{2})$.

\end{lemma}

\begin{proof}

Let $k$ be as in the statement of the lemma. Then for some
neighbourhood $N $ of the diagonal in $\partial \F_{2} \times \partial \F_{2} $, $k$ is
supported on $(\partial \F_{2} \times \partial \F_{2} \times \partial \F_{2})  \backslash (N \times \partial \F_{2})$. It follows that we can
choose an extension $\tilde{k}$ of $k$ to a function on $\partial \F_{2}
\times \overline{\F}_{2} \times \overline{\F}_{2}$ such that there is a
neighbourhood $\tilde{N}$ of the diagonal in $\partial \F_{2} \times \overline{\F}_{2}$
 such
that $\tilde{k}$ is supported in $(\partial \F_{2} \times \overline{\F}_{2} \times
\overline{\F}_{2}) \backslash (\tilde{N} \times \overline{\F}_{2})$.

Now by routine
compactness arguments, it suffices to show that for $a \in \partial \F_{2}$ fixed and
$x_{n}$ a sequence in $\F_{2}$ converging to a boundary point $b \in
\partial \F_{2}$, the sequence $$\tilde{k}(x_{n}^{-1}(a ), x_{n}^{-1}, x_{n}) -
\tilde{k}(x_{n}^{-1}(a) , x_{n}^{-1} , a)$$converges to $0$ as $n \to
\infty$. We may clearly also assume without loss of generality that for all $n$ the
point $(x_{n}^{-1}(a) , x_{n}^{-1})$ lies in the complement of
$\tilde{N}$, else both terms are $0$. By Lemma 22 there exists $R>0$ such that any two points $(c,z) \in
\partial \F_{2} \times \F_{2}$ not in $\tilde{N}$ have the property that the
geodesic $[z,c)$ passes through $B_{R}(e)$. Thus, for all $n$ large enough,
$d(e, [x_{n}^{-1} , x_{n}^{-1}(a)) \le R$. But then $d(x_{n} , [e,a))
\le R$ for all $n$. If a sequence in a hyperbolic space
remains at fixed, bounded distance from a geodesic ray, it must
converge to the endpoint of the ray. Hence $x_{n}
\to a$, and we are done by continuity of $\tilde{k}$ in the third
variable.

\end{proof}

\begin{proof} (Of Proposition 24). Consider the class $[i\otimes 1_{A}]\otimes_{A\otimes A\otimes A} \bigl( 1_{A} \otimes
    \sigma_{12}^{*}\Delta\bigr)$, which is represented by the map $B
  \otimes A \to \mathcal{Q}(A \otimes l^{2}\F_{2})$ in Equation 6.

 Define a unitary map of Hilbert
modules $U: A \otimes l^{2}\F_{2} \to A\otimes l^{2}\F_{2}$ by the formula $U (a \otimes
e_{x}) = x \cdot a \otimes e_{x}$. Let $\mathrm{Ad}_{U}$ denote
  the inner automorphism of $\mathcal{Q}(A \otimes l^{2}\F_{2})$ given by $\pi(T)
  \mapsto\pi ( UTU^{*})$ and let $\iota '$ denote the homomorphism $B \otimes A \to \mathcal{Q}(A
  \otimes l^{2}\F_{2})$ $$\iota ' (a_{1} \otimes a_{2} \otimes a_{3}) =  \mathrm{Ad}_{U}\bigl(a_{1}\otimes  \rho(a_{2}) \lambda (a_{3})\bigr).$$  We claim that
  $\iota ' = \iota.$ It is a simple matter to check that $\iota_{|_{B \otimes C_{r}^{*}(\F_{2})}} =
  \iota ' _{|_{B \otimes C_{r}^{*}(\F_{2})}}$, where $B \otimes
  C_{r}^{*}(\F_{2})$ is viewed as a sub-algebra of $B \otimes A$, and
  that for $b \in B$ and $f \in C(\partial \F_{2})$, we have $\iota (b \otimes f) = \tau (b) \pi (f \otimes 1)$ whereas
  $\iota ' (b \otimes f) = \tau (b)\bigl( 1 \otimes \lambda (f)\bigr)$. Thus it remains to prove that $\tau (b) \pi \bigl( 1\otimes  M_{\tilde{f}} -
  f\otimes 1\bigr) = 0$ in the Calkin algebra $\mathcal{Q}(A \otimes l^{2}\F_{2})$
  whenever $b \in B$, $f \in C(\partial \F_{2})$ and $\tilde{f}$ is an
  extension of $f$ to $\overline{\F}_{2}$. The collection of $b$ of
  the form $\sum \gamma F_{\gamma}$ with each $F_{\gamma} \in
  C_{c}(\partial^{2}\F_{2})$ is dense in $B$, and hence it suffices to
  prove the result for $b$ having this form. Hence it is sufficient to prove the
  result for $b = F \in C_{c}(\partial^{2}\F_{2})$. We are now done by
  Lemma 25 with $k(a,b,c) = F(a,b)f(c)$.

\end{proof}

\section{conclusion of the proof}

Now consider the class $[i \otimes 1_{A}]\otimes_{A\otimes A \otimes
  A} (1_{A} \otimes
\sigma_{12}^{*}\Delta)$, which we have shown has the form $[\iota]$,
where $\iota$ is as in Proposition 24. We are interested in calculating the
Kasparov product of the class of this extension, and the class $[D]
\otimes 1_{A} \in KK^{-1}(A , B \otimes A)$.

Recall that $$[D] = [\hat{d}_{\R}] \otimes_{C_{0}(\R)} [v-\chi]$$where
$[v-\chi]$ is the class of the homomorphism $C_{0}(\R) \to B$ induced
by mapping $\psi$ to $v-\chi$. 

Hence $[D] \otimes 1_{A} = ([\hat{d}_{\R}] \otimes 1_{A})
\otimes_{C_{0}(\R) \otimes A} ([v-\chi] \otimes 1_{A})$, where
$[v-\chi] \otimes 1_{A}$ is represented by the homomorphism $C_{0}(\R)
\otimes A \to B \otimes A$ induced by mapping $\psi \otimes a \mapsto
(v-\chi) \otimes a$. 

The Kasparov product $$([D]\otimes 1_{A}) \otimes_{B\otimes A}
[i\otimes 1_{A}]\otimes_{A\otimes A\otimes A} (1_{A} \otimes \sigma_{12}^{*}(\Delta)$$ therefore
has the form $$([\hat{d}_{\R}] \otimes 1_{A}) \otimes_{C_{0}(\R)
  \otimes A} \bigl( ([v-\chi] \otimes 1_{A}) \otimes_{B\otimes
  A}[\iota]\bigr)$$and $([v-\chi]\otimes 1_{A}) \otimes_{B\otimes A}
[\iota]$ is represented by the homomorphism $C_{0}(\R) \otimes A \to
\Q (A \otimes l^{2}\F_{2})$ induced by mapping $$\psi \otimes a \mapsto
\tau (v-\chi) \pi (\varphi (a)).$$But this homomorphism has the form
stated in the hypothesis of Lemma 8. By that lemma, $$([\hat{d}_{\R}]
\otimes 1_{A}) \otimes_{C_{0}(\R) \otimes A} ([v-\chi]\otimes 1_{A})
\otimes_{B\otimes A} [\iota]$$ is represented by the $KK(A,A)$ cycle
$\bigl( A \otimes l^{2}\F_{2} , \bar{U} +1\bigr)$, where $\bar{U}$ is any lift to
$\B (A \otimes l^{2}\F_{2})$ of the element $\tau (v-\chi) \in \Q(A
\otimes l^{2}\F_{2})$; where the Hilbert $(A,A)$-bimodule $A \otimes l^{2}\F_{2}$ has its
standard right $A$-module structure; and where it has the left $A$-module structure
given by the homomorphism $\varphi: A \to \B (A \otimes
l^{2}\F_{2})$. 

In particular, the {\it bimodule} is in fact the same as the bimodule appearing in the Julg and Vallette cycle appearing in
Lemma
21. 

It remains to calculate a lift $\bar{U}$ of $\tau (v-\chi)$ and show
that in fact such a lift can be chosen which agrees with the operator
$\bar{W}$ of Lemma 21.

We first construct a lift of $\tau (v)$. Recall that $v = \sum_{\gamma
  \in S} \chi_{\gamma}\gamma \otimes (1-\chi_{\gamma})\gamma$, where
  $S$ is a basis for $\F_{2}$. Each $\gamma$ is mapped under $\tau$ to
  the image in the Calkin algebra of the right translation operators
  $1\otimes \rho (\gamma): A \otimes l^{2}\F_{2} \to A \otimes
  l^{2}\F_{2}$. Consider each term $F_{\gamma} = \chi_{\gamma} \otimes
  (1-\chi_{\gamma}) \in C_{c}(\partial^{2}\F_{2})$. Let
  $\tilde{\chi}_{\gamma}$ denote the function on $\F_{2}$ given by 

$$\tilde{\chi}_{\gamma}(g) =
\begin{cases}  1 &  \quad  \gamma \in [e,g]  \\0& \quad \mathrm{else}
\end{cases}.$$

Then $\tilde{\chi}_{\gamma}$ extends continuously to
$\overline{\F}_{2}$ and the restriction of $\tilde{\chi}_{\gamma}$ to
$\partial \F_{2}$ is $\chi_{\gamma}$. Let then $$\tilde{F}_{\gamma} =
\chi_{\gamma} \otimes (1-\tilde{\chi}_{\gamma}),$$ which is an
extension to $\partial \F_{2} \times \overline{\F}_{2}$ of
$F_{\gamma}$. Forming $\tilde{F}_{\gamma}'$ as per the recipe
described in the definition of $\tau$, we obtain the function
$$\tilde{F}_{\gamma}'(a,g) = \tilde{F}_{\gamma}(g^{-1}a , g^{-1})  =
\begin{cases}  1 &  \quad \gamma \in [e,g^{-1}a)\; \mathrm{ and}\; \gamma
  \notin [e, g^{-1}] \\0& \quad \mathrm{else} \end{cases}.$$

We remind the reader that the statement: ``$x \in [e,y]$,'' for $x, y \in
\F_{2}$ may be equivalently read: ``the reduced expression of $y$
contains $x$ as an initial subword,'' or more shortly, ``$y$ begins
with $x$.''

With this in mind, consider the first case above. If $g^{-1}a$ begins with $\gamma$ but $g^{-1}$ does
not, it follows there is cancellation between $g^{-1}$ and $a$; more
precisely, $a$ must begin with $g$, followed by $\gamma$. (Since
$g^{-1}$ does not begin with $\gamma$, $g$ does not end in
$\gamma^{-1}$, and hence $g\gamma$ is in fact reduced.) We have:

$$\tilde{F}_{\gamma}'(a,g) =
\begin{cases}  1 &  \quad   g\gamma \in [e,a), \; \mathrm{and}\; g\;  \mathrm{does
    \; not \; end \; in}\;
  \gamma^{-1} \\0& \quad \mathrm{else} \end{cases}.$$ 

Now consider the {\it operator} $\tilde{F}_{\gamma}' \, (1 \otimes
v_{\gamma}) \in A \otimes \B (l^{2}\F_{2}) \subset \B (A \otimes
l^{2}\F_{2})$. This operates by $$(\sum f_{h}h)\otimes e_{g} \mapsto
(\sum \tilde{F}_{\gamma}'(\; \cdot \; , g\gamma^{-1})f_{h}h)\otimes
e_{g\gamma^{-1}}$$for $\sum f_{h}h$ an arbitrary element of the cross
product $A$. From our above work, we see that $\tilde{F}_{\gamma}'(\;
\cdot  \; , g\gamma^{-1}) = 0$ unless $g$ ends in $\gamma$. On the
other hand, if $g$ does end in $\gamma$, $g \gamma^{-1}$ does not end
in $\gamma^{-1}$. Hence we see that the above operator sends $$(\sum
f_{h}h) \otimes e_{g} \mapsto  \begin{cases}   (\sum
\chi_{g}  f_{h}h) \otimes e_{g\gamma^{-1}}  &  \quad   g \; \mathrm{ends \; in \; }\gamma \\0& \quad \mathrm{else} \end{cases}.$$ 

We see
finally, that $\bar{V} = \sum_{\gamma \in S} \tilde{F}_{\gamma}' \; \cdot \;
(1\otimes v_{\gamma})$, which is a lift of $\tau (v)$, acts on $A
\otimes l^{2}\F_{2}$ by $$(\sum f_{h}h) \otimes e_{g} \mapsto (\sum
\chi_{g} f_{h}h) \otimes e_{g'},$$ where the prime notation is as in the
discussion just prior to Lemma 20.

In particular, $\bar{V}$ as an operator on $A \otimes l^{2}\F_{2}$, where
the latter is regarded as functions $\F_{2} \to C(\partial \F_{2})
\otimes l^{2}\F_{2}$, has the form $$(\bar{V}\xi)(g) = V ( \xi
(g)),$$where $V$ is the operator $C(\partial \F_{2}) \otimes
l^{2}\F_{2} \to C(\partial \F_{2}) \otimes
l^{2}\F_{2}$, $$V(f \otimes e_{g}) = \chi_{g} f \otimes
e_{g'}.$$Otherwise expressed, let $\xi$ be an element of $C(\partial \F_{2} ;
l^{2}\F_{2})$ of the form $\xi (a) = \sum \xi_{g}(a) e_{g}$, where
each $\xi_{g}$ is a scalar-valued function on $\partial \F_{2}$. Then
$$(V\xi )(a) = \sum_{g \in [e,a)} \xi_{g}(a) \otimes e_{g'}.$$

Now apply the same calculations to the element $\tau (\chi)$. We
obtain the operator (projection) $\bar{P}$ on $A \otimes l^{2}\F_{2}$
given by $\bar{P} = \sum \tilde{F}_{\gamma}' \in A \otimes \B
(l^{2}\F_{2}) \subset \B (A \otimes l^{2}\F_{2})$. 

We have: $\bar{V} - \bar{P}$ is an operator whose projection to the
Calkin algebra is $\tau (v-\chi)$ as required. Let it be denoted by
$\bar{U}$. Form $\bar{F} = \bar{U}+1$.

Our calculations show that $\bar{F}$ is an operator having the form
$$(\bar{F}\xi) (g) = F(\xi (g)),$$ where $F: C(\partial \F_{2})
\otimes l^{2}\F_{2} \to  C(\partial \F_{2})
\otimes l^{2}\F_{2}$ is the operator $$(F\xi )(a) = \sum_{g \in [e,a)}
\xi_{g}(a) \otimes e_{g'} + \sum_{g \notin [e,a)} \chi_{g}(a) \otimes
e_{g},$$ which is precisely the operator $W$ of Lemma 20. That is, $F
= W$ and therefore $\bar{F} = \bar{W} \in \B(A \otimes l^{2}\F_{2})$. 

We are now done, having shown by direct computation that$$([D] \otimes 1_{A} ) \otimes_{B\otimes A} [i\otimes
1_{A}] \otimes_{A \otimes A \otimes A} (1_{A} \otimes
\sigma_{12}^{*}\Delta) = [\bigl( C(\partial \F_{2}) \rtimes \F_{2} \otimes
l^{2}\F_{2} , \bar{W} \bigr)] = \lambda (\gamma_{\partial \F_{2} \rtimes \F_{2}})
$$and therefore that $$([D] \otimes 1_{A} ) \otimes_{B\otimes A} [i\otimes
1_{A}] \otimes_{A \otimes A \otimes A} (1_{A} \otimes
\sigma_{12}^{*}\Delta) = 1_{A}.$$

\end{document}